\documentclass[11pt]{amsart}
\usepackage{amsfonts}
\usepackage{mathrsfs}
\begin{document}
\baselineskip=20pt

\title{Classification of Graded Left-symmetric Algebra Structures on Witt
and Virasoro Algebras}

\author{Xiaoli Kong}
\address{Xiaoli Kong\\
School of Mathematical Sciences, Xiamen University,
 Xiamen, Fujian 361005, P.R.
China}\email{kongxl.math@gmail.com}

\author{Hongjia Chen}
\address{Hongjia Chen\\
Department of Mathematics, University of Science and Technology of China,
 Anhui, Hefei 230026, P.R. China}
\email{hjchen@mail.ustc.edu.cn}

\author{Chengming Bai}
\address{Chengming Bai\\
Chern Institute of Mathematics and LPMC, Nankai University, Tianjin 300071, P.R.China}
\email{baicm@nankai.edu.cn}

\def\shorttitle{Left-symmetric algebra structures on Witt and
Virasoro algebras}

\begin{abstract} We find that a compatible graded left-symmetric algebra
structure on the Witt algebra induces an indecomposable module of
the Witt algebra with 1-dimensional weight spaces by its left
multiplication operators. From the classification of such modules of
the Witt algebra, the compatible graded left-symmetric algebra
structures on the Witt algebra are classified. All of them are
simple and they include the examples given by Chapoton and
Kupershmidt. Furthermore, we classify the central extensions of
these graded left-symmetric algebras which give the compatible
graded left-symmetric algebra structures on the Virasoro algebra.
They coincide with the examples given by Kupershmidt.
\end{abstract}


\maketitle

\section{Introduction}
Left-symmetric algebras (or under other names like pre-Lie algebras,
quasi-associative algebras, Vinberg algebras and so on) are
Lie-admissible algebras (see Proposition 2.2). They were introduced
by A. Cayley in 1896 as a kind of rooted tree algebras (\cite{Ca}).
They also arose from the study of convex homogenous cones
(\cite{V}), affine manifolds and affine structures on Lie groups
(\cite{Ko}), deformation of associative algebras (\cite{G}) in
1960s. As it was pointed out in \cite{CL}, the left-symmetric
algebra ``deserves more attention than it has been given''. They
appear in many fields of mathematics and mathematical physics. In
\cite{Bu2}, there is a survey of certain different fields where
left-symmetric algebras play an important role, such as vector
fields, rooted tree algebras, words in two letters, vertex algebras,
operad theory, deformation complexes of algebras, convex homogeneous
cones, affine manifolds, left-invariant affine structures on Lie
groups. In addition, left-symmetric algebras have close relations
with symplectic and complex structures on Lie groups and Lie
algebras (\cite{Chu}, \cite{LM}, [DaM1-2], \cite{AS}), phase spaces
of Lie algebras ([Ku1-2], \cite{Ba}), certain integrable systems
(\cite{Bo}, \cite{SS}), classical and quantum Yang-Baxter equations
(\cite{Ku3}, \cite{ES}, \cite{GS}, \cite{DiM}), combinatorics
(\cite{E}) and so on. In particular, they play a crucial role in the
Hopf algebraic approach of Connes and Kreimer to renormalization of
perturbative quantum field theory (\cite{CK}).

It is not easy to study left-symmetric algebras. Due to the
nonassociativity, there is neither a suitable representation theory
nor a complete structure theory like other classical algebras such
as associative algebras and Lie algebras. For example, it is far
from the classification of semisimple left-symmetric algebras (a
left-symmetric algebra is called simple if it has not nontrivial
ideals and a semisimple left-symmetric algebra is a direct sum of
simple ones). In fact, the classification of complex simple
left-symmetric algebras are only in low dimensions and some special
cases in certain higher dimensions (\cite{Bu1}).

On the other hand, an important approach to study left-symmetric
algebras is through the representation theory of their sub-adjacent
Lie algebras. It is known that there exists a compatible
left-symmetric structure on a Lie algebra ${\mathcal G}$ if and only
if ${\mathcal G}$ has an \'etale affine representation or
equivalently, ${\mathcal G}$ has a bijective 1-cocycle associated to
a representation (cf. \cite{Me}, \cite{Ki}). Unfortunately, the
sub-adjacent Lie algebra of a finite-dimensional left-symmetric
algebra over a field of characteristic zero cannot be semisimple. So
the beautiful representation theory of finite-dimensional semisimple
Lie algebras cannot be used here. However, in the case of the
infinite dimensional, there exists a semisimple Lie algebra with a
compatible left-symmetric algebra structure. One of such examples is
given as a left-symmetric Witt algebra, which is regarded as the
first important example of infinite-dimensional left-symmetric
algebras (it was also regarded as one of the origins of
left-symmetric algebras in \cite{Bu2}). Let $U$ be an associative
commutative algebra, and ${\mathcal D}=\{
\partial_1,\cdots,\partial_n\}$ be a system of commutating
derivations of $U$. Then the vector space
$${\rm Vec}(n)=\{ \sum_{i=1}^n u_i\partial_i\mid u_i\in U, \partial_i\in
{\mathcal D}\},\eqno (1.1)$$ is called a left-symmetric Witt algebra
under the multiplication ([Ku1], \cite{DL}, \cite{Bu2})
$$u\partial_i\circ v\partial_j=u\partial_i(v)\partial_j.\eqno (1.2)$$
In particular, when $n=1$, $U={\mathbb{F}}[x,x^{-1}]$ and
$\partial=\frac{\partial}{\partial x}$, the $W^1={\rm Vec}(1)$ is a
simple left-symmetric algebra, whose sub-adjacent Lie algebra is the
Witt algebra.

Therefore, it is natural to consider the classification of the
compatible left-symmetric algebra structures on the Witt algebra. It
is easy to know that all compatible left-symmetric algebras on the
Witt algebra are simple. Since the Witt algebra is graded, it is
also natural to suppose that the compatible left-symmetric algebras
should be graded. Hence, in this paper, we mainly consider the
left-symmetric algebras with a basis $\{x_n\mid n\in {\mathbb Z}\}$
satisfying
$$x_ix_j=f(i,j)x_{i+j},\;\;[x_i,x_j]=x_ix_j-x_jx_i=(j-i)x_{i+j},\eqno (1.3)$$
where $f$ is a complex-value function on ${\mathbb Z}\times {\mathbb
Z}$. In particular, $W^1$ is just the case $f(i,j)=1+j$ and it is a
Novikov algebra, which is a left-symmetric algebra with commutative
right multiplication operators (\cite{GD}, \cite{BN}).

Moreover, like $W^1$, some other special cases satisfying equation
(1.3) have been already discussed. For example, Chapoton in
\cite{Cha} gave the classification in the case $f(i,j)=g(i)h(j)$ (up
to a change of basis). As it was pointed out in \cite{Cha}, ``this
Ansatz for the product has no special meaning, except that it allows
for a full classification". Kupershmidt gave a solution
$f(i,j)=\frac{j(1+bj)}{1+b(i+j)}$ for $b=0$ or {$b^{-1}\notin
{\mathbb Z}$} in \cite{Ku2}. Osborn classified a class of simple
infinite dimensional Novikov algebras which includes a
classification of Novikov algebras satisfying equation (1.3), that
is, $f(i,j)=\alpha+j$ for any {$\rm{Re}\alpha> 0$ or
$\rm{Re}\alpha=0,\rm{Im}\alpha \geq 0$} (\cite{O}). However, to our
knowledge, a complete classification of left-symmetric algebras
satisfying equation (1.3) is still unknown.

In this paper, we will mainly use the representation theory of the
Witt algebra. We find that a regular representation of the Witt
algebra induced by a compatible left-symmetric algebra structure
satisfying equation (1.3) (through its left multiplication
operators) belongs to a class of important modules of the Witt
algebra, which have been classified. Hence such left-symmetric
algebras can be classified through studying the relations between
them. As it was done in Lie algebras, it is natural to consider the
central extensions of the left-symmetric algebras satisfying
equation (1.3), by which we can get the compatible left-symmetric
algebra structures on the Virasoro algebra. It is interesting to see
that there exist only one class of non-trivial central extensions
which coincides with a result in \cite{Ku2}.

We would like to point out some interesting remarks from the
following aspects.

$\circ$\quad Our methods in this paper provide an approach to the
study of the possible compatible left-symmetric algebra structures
on the infinite-dimensional (semi)simple Lie algebras with a good
representation theory.

$\circ$\quad The left-symmetric algebras obtained in this paper are
the simple graded left-symmetric algebras of growth one. It is
difficult to classify all of it. Our study can be regarded as the
first step, as Mathieu solved the analogous problem for Lie algebras
([Ma1-3]).

$\circ$\quad It will be also interesting to consider their relations
with the vertex (operator) algebras which are the fundamental
algebraic structures in conformal field theory, since a vertex
algebra is equivalent to a left-symmetric algebra and a Lie
conformal algebra with some compatible conditions (\cite{BK},
\cite{LL}).

This paper is organized as follows. In Section 2, we give some
necessary definitions, notations and basic results on left-symmetric
algebras and the representation theory of the Witt and Virasoro
algebras. In Section 3, we prove that a left-symmetric algebra
structure satisfying equation (1.3) induces an indecomposable module
of the Witt algebra with 1-dimensional weight spaces by its left
multiplication operators. Therefore, such left-symmetric algebras
are classified through the classification of those modules. In
Section 4, the non-trivial central extensions of the left-symmetric
algebras obtained in Section 3 are discussed. A classification of
the compatible left-symmetric algebra structures on the Virasoro
algebras are obtained.

Throughout this paper, all algebras are over the complex field
${\mathbb C}$ and the indices $m, n, l, i, j, k \in\mathbb{Z}$,
unless otherwise stated.

\section{Preliminaries and fundamental results}

{\bf Definition 2.1.} Let $A$ be a vector space over a field
$\mathbb{F}$ equipped with a bilinear product $(x,y)\rightarrow xy$.
$A$ is called a left-symmetric algebra if for any $x,y,z\in A$, the
associator
$$(x,y,z)=(xy)z-x(yz)\eqno (2.1)$$
is symmetric in $x,y$, that is,
$$(x,y,z)=(y,x,z),\;\;{\rm or}\;\;{\rm
equivalently}\;\;(xy)z-x(yz)=(yx)z-y(xz).\eqno (2.2)$$

Left-symmetric algebras are Lie-admissible algebras (cf. [Me]).

{\bf Proposition 2.2.} Let $A$ be a left-symmetric algebra. For any
$x\in A$, let $L_x$ denote the left multiplication operator, that
is, $L_x(y)=xy$ for any $y\in A$. Then we have the following
results:

(1) The commutator
$$[x,y]=xy-yx,\;\;\forall x,y\in A,\eqno (2.3)$$
defines a Lie algebra ${\mathcal G}(A)$, which is called a
sub-adjacent Lie algebra of $A$ and $A$ is also called a compatible
left-symmetric algebra structure on the Lie algebra ${\mathcal
G}(A)$.

(2) Let $L:{\mathcal G}(A)\rightarrow gl(A)$ with $x\mapsto L_x$.
Then $(L,A)$ gives a  representation of the Lie algebra ${\mathcal
G}(A)$, that is,
$$[L_x,L_y]=L_{[x,y]},\;\;\forall x,y\in A. \eqno (2.4)$$
We call it a regular representation of the Lie algebra ${\mathcal
G}(A)$.

There is not a compatible left-symmetric algebra structure on any
Lie algebra ${\mathcal G}$. A sufficient and necessary condition for
a Lie algebra with a compatible left-symmetric algebra structure is
given as follows. Let ${\mathcal G}$ be a Lie algebra and
$\rho:{\mathcal G}\rightarrow gl(V)$ be a representation of
${\mathcal G}$. A 1-cocycle $q:{\mathcal G}\rightarrow V$, the
linear map on vector space associated to $\rho$  (denoted by
$(\rho,q)$) satisfying
$$q[x,y]=\rho(x)q(y)-\rho(y)q(x),\ \forall x,y\in {\mathcal G}.\eqno (2.5)$$
Let $A$ be a left-symmetric algebra and $\rho: {\mathcal
G}(A)\rightarrow gl(V)$ be a representation of its sub-adjacent Lie
algebra. If $g$ is a homomorphism of the representations from $A$ to
$V$, then $g$ is a 1-cocycle of ${\mathcal G}(A)$ associated to
$\rho$ .

{\bf Proposition 2.3.} Let ${\mathcal G}$ be a Lie algebra. Then
there is a compatible left-symmetric algebra structure on ${\mathcal
G}$ {if and only if} there exists a bijective 1-cocycle of
${\mathcal G}$.

In fact, let $(\rho,q)$ be a bijective 1-cocycle of ${\mathcal G}$,
then
$$x*y=q^{-1}\rho(x)q(y),\;\;\forall x,y\in {\mathcal G},\eqno (2.6) $$
defines a left-symmetric algebra structure on ${\mathcal G}$.
Conversely, for a left-symmetric algebra $A$, the identity
transformation $id$ is a 1-cocycle of ${\mathcal G}(A)$ associated
to the regular representation $L$.

On the other hand, we recall the definition of the Witt algebra $W$
and Virasoro algebra $\mathscr{V}$. The Witt algebra $W$ (of rank
one) is a complex Lie algebra with a basis $\{x_n\mid n\in
\mathbb{Z}\}$ whose commutation relations satisfy
$$[x_m,x_n]=(n-m)x_{m+n}.\eqno (2.7)$$
 The Virasoro algebra $\mathscr{V}$ is the
central extension of $W$ defined by the $2$-cocycle
$$\Omega(x_m,x_n) =\frac{1}{12}(n^3-n)\delta_{m+n,0}.$$ That is,
$\mathscr{V}$ is a complex Lie algebra with a basis $\{\theta,
x_n\mid n\in \mathbb{Z}\}$ whose commutation relations satisfy
$$[\theta,
x_n]=0,\;\;[x_m,x_n]=(n-m)x_{m+n}+\delta_{m+n,0}\frac{n^3-n}{12}\theta.\eqno
(2.8)$$

Since an ideal of a left-symmetric algebra is still an ideal of its
sub-adjacent Lie algebra and $W$ is a simple Lie algebra, we have
the following conclusion.

{\bf Proposition 2.4.} Any compatible left-symmetric algebra
structure on the Witt algebra $W$ is simple.

For any module $V$ of a Lie algebra ${\mathcal G}$, the action of
${\mathcal G}$ on $V$ is denoted by $xv$ for any $x\in {\mathcal G}$
and $v\in V$. Recall that a module $V$ of a Lie algebra is called
indecomposable if $V$ cannot be decomposed into a direct sum of two
proper submodules. A weight space of a module $V$ of the Witt
algebra $W$ or the Virasoro algebra $\mathscr{V}$ is the non-zero
vector space $V_\lambda$ defined by
$$V_\lambda=\{x\in V\mid x_0v=\lambda v\}.\eqno (2.9)$$

For later use, we give the description of the indecomposable modules
with 1-dimensional weight spaces of the Virasoro algebra
$\mathscr{V}$ (\cite{KS}, [MP1-2], \cite{Ma4}).

$-$The $\mathscr{V}$-module $A_{\alpha,\beta}$ of Feigin-Fuchs with
$\alpha,\beta\in\mathbb{C}$ and $0\leq \rm{Re} \alpha<1$ , whose
action on a basis $\{v_n\mid n \in \mathbb{Z} \}$ is given by
$$x_iv_n=(\alpha+n+i\beta)v_{n+i},\quad \theta v_n=0.$$

$-$The maximal proper  $\mathscr{V}$-submodule of $A_{0,1}$, called
$A'_{0,1}=A_{0,1}\setminus \mathbb{C}v_0\;\ (A_{0,1}/A'_{0,1}$ is
trivial and $A_{0,0}/\mathbb{C}v_0 \simeq A'_{0,1})$,\ \ whose
action on a basis $\{v_n\mid\ n\in \mathbb{Z}\setminus\{0\}\}$ is
given by
$$x_iv_n=(n+i)v_{n+i},\quad \theta v_n=0,\quad \forall \  n\neq0.$$

$-$ The $\mathscr{V}$-module $A_\alpha$ with $\alpha\in \mathbb{C}$
whose action on a basis $\{v_n\mid n \in \mathbb{Z} \}$ is given  by
\begin{eqnarray*}
&&x_iv_n=(n+i)v_{n+i}, \quad\forall \  n\neq0,\\
&&x_iv_0=i(\alpha+i)v_i,\quad \theta v_n=0.
\end{eqnarray*}

$-$ The $\mathscr{V}$-module $B_\beta$ with $\beta\in \mathbb{C}$
whose action on a basis $\{v_n\mid n \in \mathbb{Z} \}$ is given  by
\begin{eqnarray*}
&&x_iv_n=nv_{n+i}, \quad \mbox{if  }n+i\neq0,\\
&&x_iv_{-i}=-i(\beta+i)v_0,\quad \theta v_n=0.
\end{eqnarray*}

{\bf Theorem 2.5.}(\cite{KS}) Any indecomposable nontrivial module
of the Virasoro algebra $\mathscr{V}$ with 1-dimensional weight
spaces is isomorphic to one of $A_{\alpha,\beta}, A'_{0,1},
A_\alpha,$  $B_\beta$.

Obviously, $A_{\alpha,\beta}, A'_{0,1}, A_\alpha,$ and $B_\beta$ are
also the modules of the Witt algebra $W$.

{\bf Corollary 2.6.} Any indecomposable nontrivial module of the
Witt algebra $W$ with 1-dimensional weight spaces is isomorphic to
one of $A_{\alpha,\beta}, A'_{0,1},$ $A_\alpha,$  $B_\beta$.

\section{Compatible left-symmetric algebra structures on the Witt algebra}

As we said in the introduction, we study the following compatible
left-symmetric algebra structure on the Witt algebra $W$ with the
multiplication
$$x_mx_n=f(m,n)x_{m+n}, \eqno(3.1)$$
where $\{ x_n\mid n\in {\mathbb Z}\}$ is a basis of $W$ satisfying
equation (2.7). To avoid confusion, we denote such a left-symmetric
algebra structure by $V$. Then $V$ is a regular module of $W$
defined by the left multiplication operators of $V$, that is, in
equation (3.1), we let $x_m\in W$ and $x_n, x_{m+n}\in V$.

In fact, $V$ is a compatible left-symmetric algebra structure on $W$
if and only if
$$[x_m,x_n]=x_mx_n-x_nx_m=(n-m)x_{m+n}\mbox{ and }\
(x_m,x_n,x_l)=(x_n,x_m,x_l).\eqno(3.2)$$ They hold if and only if
$f(m,n)$ satisfies the following equations:
$$f(m,n)-f(n,m)=n-m,\eqno(3.3)$$
$$(n-m)f(m+n,l)=f(n,l)f(m,n+l)-f(m,l)f(n,m+l).\eqno(3.4)$$

{\bf Lemma 3.1.} $f(m,0)=f(0,0)$.

{\bf Proof.} Let $n=l=0$ in equation (3.4), we have
$$-mf(m,0)=f(0,0)f(m,0)-f(m,0)f(0,m).$$ Then by equation (3.3), we
have
$$f(m,0)(f(0,m)-m-f(0,0))=f(m,0)(f(m,0)-f(0,0))=0.$$
Therefore, $f(m,0)=0\ $ or $\ f(m,0)=f(0,0).$

If $f(0,0)=0$, then $f(m,0)=0=f(0,0)$. The conclusion holds.

If $f(0,0)\neq 0$, set
$$I_1=\{m\in {\mathbb Z} \mid f(m,0)=0\},\ I_2=\{m\in {\mathbb Z}\mid f(m,0)=f(0,0)\}.$$
Obviously, $I_1\cup I_2={\mathbb Z}$, $I_1\cap I_2=\emptyset$ and
$0\in I_2$.

Let $l=0$ in equation (3.4), we have
$$(n-m)f(m+n,0)=f(n,0)f(m,n)-f(m,0)f(n,m).\eqno(*)$$
Then we obtain the following results:

 (a) If  $m,n\in I_s, m\neq n$ $(s=1,2)$, then $m+n\in I_s$;

 (b)  If $m\in I_1$, then $-m\in I_2$.

\noindent Next we prove that $1,-1,2,-2\notin I_1$, that is
$1,-1,2,-2\in I_2$. Thus by the above result (a), it is easy to know
that $I_2=\mathbb{Z}$ and $I_1=\emptyset$. The conclusion holds.

{\it  (I).} $1\notin I_1$.

 Otherwise, we suppose $1\in I_1$. Hence $-1\in I_2$, $2\in I_1$ (if
$2\in I_2$, then $1=2+(-1)\in I_2$ which is a contradiction).
Therefore, for any $n>0$, we know $n\in I_1$ by induction. Thus
$$I_1=\{1,2,\ldots\},\,\;\;
I_2=\{0,-1,-2,\ldots\}.$$

 If $-n \geq m>0$, then equation ($*$)
becomes $$(n-m)f(0,0)=f(0,0)f(m,n).$$ Hence $f(m,n)=n-m,\ f(n,m)=0$.
In particular, if  $n<0$, we have $f(n,1)=0$.

If $m>-n>0$, then $m+n>0$, and equation ($*$) becomes
$$f(0,0)f(m,n)=0.$$
Hence $f(m,n)=0,\ f(n,m)=m-n$. In particular, if  $m\geq -n>0$,
$f(n,m+1)=m+1-n$.

Let $l=1$ in equation (3.4), we have
$$(n-m)f(m+n,1)=f(n,1)f(m,n+1)-f(m,1)f(n,m+1).$$
Therefore when $m\geq -n>0$, the above equation becomes
$$(n-m)f(m+n,1)=(n-m-1)f(m,1).$$
Let $m=-n>0$. Then we have
$$2mf(0,1)=(2m+1)f(m,1),\;\;\; m \geq 1,$$
where $f(0,1)=1$ since $1\in I_1$ and $f(1,0)=0$.

\noindent Let $m=1-n>1$. Then we have
$$(2m-1)f(1,1)=2mf(m,1),\;\; \ m \geq
2.$$ Thus, by the above two equations, we know that
$$f(1,1)= \frac{4m^2}{4m^2-1}$$ for all $m \geq 2$, which
is impossible.

{\it (II).} $-1\notin I_1$. The proof is similar to the proof in
Case (I) by symmetry.

{\it (III).} $2\notin I_1$.

Otherwise suppose $2\in I_1$. Therefore $1,0,-1,-2\in I_2$. Now for
any $n\geq 2$, we have $n\in I_1$ by induction since $n=(n+1)+(-1)$.
Thus
$$I_1=\{2,3,\ldots\},\ I_2=\{1,0,-1,-2,\ldots\}.$$

If $m-1>-n>0$, that is, $m+n>1,\ m>1$, then equation ($*$) becomes
$$0=f(0,0)f(m,n).$$
Hence $f(m,n)=0,\ f(n,m)=m-n$. In particular, if $m\geq -n\geq 2$,
then we have $f(m+2,n)=0,\ f(n,m+2)=m+2-n$.

If $-m\geq n>1$, then equation ($*$) becomes
$$(n-m)f(0,0)=-f(0,0)f(n,m).$$
Hence $f(n,m)=m-n,\ f(m,n)=0$. In particular, if $n\leq -2$, then
$f(n,2)=0$.

Let $l=2$ in equation (3.4), we have
$$(n-m)f(m+n,2)=f(n,2)f(m,n+2)-f(m,2)f(n,m+2).$$
Therefore, when $m\geq -n\geq 2$, the above equation becomes
$$(m-n)f(m+n,2)=(m+2-n)f(m,2).$$
Let $m=-n \geq 2$. Then we have
$$2mf(0,2)=(2m+2)f(m,2),\;\; m \geq 2,$$
where $f(0,2)=2$ since $2\in I_1$ and $f(2,0)=0$.

\noindent Let $m=1-n \geq 3$. Then we have
$$(2m-1)f(1,2)=(2m+1)f(m,2), \;\; m \geq 3.$$
Thus, by the above two equations, we know that
$$f(1,2)= \frac{4m^2+2m}{2m^2+m-1}$$
for all $m \geq 3$, which is impossible.

{\it Case (IV).} $-2\notin I_1$. The proof is similar to the proof
in Case (III) by symmetry. \hfill$\Box$

{\bf Lemma 3.2.} The weight space of $V$ is $1$-dimensional.

{\bf Proof.} By Lemma 3.1 and equation (3.3), we have
$f(0,m)=f(0,0)+m$, that is, the elements $\{x_n\}$ are in different
eigenspaces of $x_0$. So every weight space of $V$ is
$1$-dimensional.\hfill$\Box$

{\bf Lemma 3.3.} $V$ is an indecomposable $W$-module.

{\bf Proof.} We assume that $V=V_1\oplus V_2$, where $V_1,V_2$ are
proper submodules of $V$.  We assume that $V_2$ is indecomposable
and $x_0\notin V_2$ without losing generality. It is well known that
any submodule of a weight module is still a weight module ([Ka]).
Then there exist two nonempty subsets $I_1,I_2$ of $\mathbb{Z}$,
such that
$$0\notin I_2,\;\; I_1\cup I_2=\mathbb{Z},\;\; I_1\cap I_2=\emptyset,\;\;{\rm
and}\;\;V_i=\bigoplus_{n\in I_i}\mathbb{C}x_n,\quad i=1,2.$$

Let $s$ be the smallest positive integer in $I_2$ (if all elements
of $I_2$ are negative, then let $s$ be the largest negative integer,
and the discussion is similar). We have the following results.

(i)\quad $f(m,0)=0, f(0,m)=m$. In fact, $x_sx_0=f(s,0)x_s\in V_1\cap
V_2$ since $x_0\in V_1$. So $f(s,0)=0$. Therefore we have
$f(m,0)=f(s,0)=0$ by Lemma 3.1.

(ii)\quad If $m,n\in I_2$ and $m\neq n$, then $m+n\in I_2$. In fact,
since $f(m,n)-f(n,m)=n-m\neq 0$, at least one of $f(m,n),\ f(n,m)$
is not zero. Suppose $f(m,n)\neq 0$ without losing generality. Then
$x_mx_n=f(m,n)x_{m+n}\neq 0$. Hence $x_{m+n}\in V_2$.

(iii)\quad If $m\in I_2$, then $-m\in I_1$. Otherwise, by (ii), we
know that $0=m+(-m)\in I_2$, which contradicts to our assumption.

(iv)\quad If $m\in I_k$ $(k=1,2)$, $f(n,m)\neq 0$, then $m+n\in
I_k$.

(v)\quad If $m\in I_2$ and $n\in I_1$, then $m+n\in I_1$
if and only if $f(n,m)=0$, $f(m,n)=n-m$ and $m+n\in I_2$ if and only if
$f(m,n)=0$, $f(n,m)=m-n$.

Both (iv) and (v) immediately follow from $x_nx_m=f(n,m)x_{m+n}$.

(vi)\quad $s+1\in I_2$. In fact, let $-m=n=1$ and $l=s$ in equation
(3.4), we have,
$$2f(0,s)=f(1,s)f(-1,s+1)-f(-1,s)f(1,s-1).\eqno (**)$$
By (i), $f(0,s)=s$. Since $s-1\in I_1$ due to our assumption on $s$,
$f(-1,s)=0$ by (iv). Therefore we know that $2s=f(1,s)f(-1,s+1)\neq
0$. Hence $f(1,s)\ne 0$. So $s+1\in I_2$ by (iv).

If $s=1$, then $2\in I_2$. Moreover, $\{n\in {\mathbb Z}\mid n\geq
1\} \subset I_2$ by (ii) and $\{n\in {\mathbb Z}\mid n\leq 0\}
\subset I_1$ by (iii). That is,
$$I_1=\{0,-1,-2,\ldots\},\quad I_2=\{1,2,\ldots\}.$$ By equation ($**$) and
(v), we know that $f(1,1)=\frac{2}{3}$. Let $m=-n=2$ and $l=s=1$ in
equation (3.4), we know that $f(2,1)=\frac{4}{5}$. Let
$m=2,n=-1,l=s=1$ in equation (3.4), we have
$$-3f(1,1)=f(-1,1)f(2,0)-f(2,1)f(-1,3).$$ Therefore, $2=\frac{16}{5}$,
which is a contradiction.

Now supposing that $s>1$, we have $x_{s-1}\in V_1, x_1\in V_1$.
Hence $x_1x_{s-1}\in V_1\cap V_2$ and $x_{s-1}x_1\in V_1\cap V_2$.
Therefore $f(1,s-1)=f(s-1,1)=0$. Since $f(1,s-1)-f(s-1,1)=s-2$ by
equation (3.3), we know that $s-2=0$, that is, $s=2$. Moreover,
$-1\in I_1$ by (ii), $-2\in I_1$ by (iii) and $3\in I_2$ by (vi). By
equation ($**$) and (v), we know that $f(1,2)=1$. Let $m=-n=2,
l=s=2$ in equation (3.4), we have
$$-4f(0,2)=f(-2,2)f(2,0)-f(2,2)f(-2,4).$$
By (i), $f(0,2)=2$ and $f(2,0)=0$. Hence $f(2,2)\ne 0$. Then
$4=2+2\in I_2$ by (iv). Furthermore, by (v), we know that
$f(-2,4)=6$. Therefore $f(2,2)=\frac{4}{3}$. Let $m=2,n=-1,l=s=2$ in
equation (3.4), we have
$$-3f(1,2)=f(-1,2)f(2,1)-f(2,2)f(-1,4).$$
Therefore, $3=\frac{20}{3}$, which is a contradiction, too.

Thus there does not exist such an $s$ and hence $V$ is an
indecomposable $W$-module.\hfill$\Box$

{\bf Corollary 3.4.} As a $W$-module, $V$ must be isomorphic to one
of $A_{\alpha,\beta},$ $A'_{0,1},$ $A_\alpha,$  $B_\beta$.

Next we discuss the existence and the classification of the
left-symmetric algebras $V$.

{\bf Theorem 3.5.} As a $W$-modules, $V$ is not isomorphic to
$A'_{0,1}$.

{\bf Proof.} Suppose that $V\cong A'_{0,1}$, where $V$ is given by
the multiplication $x_mx_n=f(m,n)x_{m+n}$.  Then there exists an
isomorphism $g:V\rightarrow A'_{0,1}$ of modules  such that
$g(x_0)=\sum_{i\ne 0}c_iv_i$, for a finite number of nonzero
$c_i\in\mathbb{C}$. Since $g(x_0x_0)=x_0g(x_0)$, we have
$$\sum_{i\ne 0}f(0,0)c_iv_i=\sum_{i\ne 0}c_iiv_i.$$ Because $g$ is an
isomorphism, there exists $k\neq0$, such that $c_k\neq 0$ and $
f(0,0)=k$. Then $g(x_0)=c_kv_k$ and by Lemma 3.1, $f(m,0)=k$. Since
$g(x_{-k}x_0)=x_{-k}g(x_0)$, we know that $f(-k,0)g(x_{-k})=0$.
Therefore $g(x_{-k})=0$, which is a contradiction.\hfill$\Box$

By a similar discussion as above, we have the following conclusion.

{\bf Theorem 3.6.} As a $W$-modules, $V$ is not isomorphic to
$A_\alpha$.

{\bf Theorem 3.7.} There are compatible left-symmetric algebra
structures $V_{\alpha,\epsilon}$ on the Witt algebra $W$ given by
the multiplication
$$x_mx_n= \frac{(\alpha+n+\alpha\epsilon m)(1+\epsilon
n)}{1+\epsilon(m+n)}x_{m+n},\eqno(3.5)$$ where $\alpha,\epsilon \in
\mathbb{C}$ satisfying $\epsilon=0$ or
$\epsilon^{-1}\notin\mathbb{Z}$.

{\bf Proof.} Suppose that $V \cong A_{\alpha,\beta}$, where
$\alpha,\beta\in\mathbb{C}$, $0\leq \rm{Re}\alpha<1$, and $V$ is
given by the multiplication $x_mx_n=f(m,n)x_{m+n}$. Then there
exists an isomorphism $g:V\rightarrow A_{\alpha,\beta}$ of modules
such that $g(x_0)=\sum_ic_iv_i$ for a finite number of nonzero
$c_i\in\mathbb{C}$. Since $g(x_0x_0)=x_0g(x_0)$, we have
$$\sum_if(0,0)c_iv_i=\sum_ic_i(\alpha+i)v_i.$$
So there exists $k$ such that $c_k\ne 0$ and $f(0,0)=\alpha+k$. Then
{$g(x_0)=c_kv_k$} and by Lemma 3.1, $f(m,0)=\alpha+k$.
Since $g(x_mx_0)=x_mg(x_0)$, we know that
$$f(m,0)g(x_m)=c_kx_mv_k=c_k(\alpha+k+m\beta)v_{m+k}.$$

If $\alpha+k \neq 0$, then $g(x_m)=
\frac{\alpha+k+m\beta}{\alpha+k}c_kv_{m+k}$. Moreover, in this case,
$\alpha+k+m\beta\neq0$ since $g$ is an isomorphism. Since
$g(x_mx_n)=x_mg(x_n)$, we have
$$f(m,n)g(x_{m+n})=\frac{\alpha+k+n\beta}{\alpha+k}c_kx_mv_{n+k}=
\frac{\alpha+k+n\beta}{\alpha+k}(\alpha+n+k+m\beta)c_kv_{m+n+k}.$$
Therefore,
$$f(m,n)=\frac{(\alpha+k+n+m\beta)(\alpha+k+n\beta)}{\alpha+k+(m+n)\beta}.$$
Replacing $\alpha+k$ and $ \frac{\beta}{\alpha+k}$ by $\alpha$ and
$\epsilon$ respectively, we have
$$f(m,n)= \frac{(\alpha+n+\alpha\epsilon m)(1+\epsilon
n)}{1+\epsilon(m+n)},$$ where $\alpha\neq0$ and $\epsilon=0$ or
$\epsilon^{-1} \notin \mathbb{Z}$.

If $\alpha+k=0$, then $\alpha=k=0$ since  $0\leq {\rm
Re}{\alpha}<1$.  Hence $m \beta c_kv_{m+k}=0$. Thus $\beta=0$.
Moreover, we have $f(m,0)=0$ and $f(0,m)=m$. Since
$mg(x_m)=g(x_0x_m)=x_0g(x_m)$, we know that $g(x_m)=a_m v_m$, where
$a_m \neq 0$. In particular, $a_0=c_0$. Since
$$f(m,n)a_{m+n}v_{m+n}=f(m,n)g(x_{m+n})=g(x_mx_n)=x_mg(x_n)=na_nv_{m+n},$$
we know that $f(m,n)= \frac{na_n}{a_{m+n}}$. By equation (3.3), we
have
$$(n-m)a_{m+n}=na_n-ma_m.$$
Let $m=-n$, we have $a_n+a_{-n}=2a_0.$ On the other hand, let $m=2,
n=-1$, we have $a_{-1}=3a_1-2a_2$. Therefore $a_2=2a_1-a_0$. {By
induction, we know that $a_n=na_1-(n-1)a_0$. Replacing $
\frac{a_1-a_0}{a_0}$ by $\epsilon$, we have
$$f(m,n)=\frac{na_n}{a_{m+n}}=\frac{n(na_1-na_0+a_0)}{(m+n)a_1-(m+n)a_0+a_0}=
\frac{n(1+n \epsilon)}{1+\epsilon(m+n)},$$} where $\epsilon=0$ or
$\epsilon^{-1}\notin \mathbb{Z}$. Moreover, it is just the case
$\alpha=0$ in equation (3.5).

Conversely, it is easy to know that we define a compatible
left-symmetric algebra structure on the Witt algebra $W$ by equation
(3.5).\hfill$\Box$

{\bf Theorem 3.8.} There are compatible left-symmetric algebra
structures $V^{\beta,k}$ on the Witt algebra $W$ given by the
multiplications
$$x_mx_n=(n+k)x_{m+n}, \, \mbox{if  } m+n+k \neq 0, \eqno(3.6)$$
$$x_{-n-k}x_n=\frac{(n+k)(\beta-n-k)}{\beta-k}x_{-k},  \eqno(3.7)$$
 where $\beta\in \mathbb{C}$ and $k \in \mathbb{Z}$ satisfying
$\beta -k \neq 0$.

{\bf Proof.} Suppose that $V \cong  B_ \beta$, where
$\beta\in\mathbb{C}$, and $V$ is given by the multiplication
$x_mx_n=f(m,n)x_{m+n}$. Then there exists an isomorphism
$g:V\rightarrow B_ \beta$ of modules such that $g(x_0)=\sum_ic_iv_i$
for a finite number of nonzero $c_i\in\mathbb{C}$. Since
$g(x_0x_0)=x_0g(x_0)$, we have
$$\sum_if(0,0)c_iv_i=\sum_{i\ne 0}c_iiv_i.$$
So there exists $k$ such that $c_k\ne 0$ and $f(0,0)=k$. Then
$g(x_0)={c_kv_k}$ and by Lemma 3.1, $f(m,0)=k$.

{\it Case (I).} $k\neq0$. Since $g(x_mx_0)=x_mg(x_0)$, we have
$$f(m,0)g(x_m)=c_kkv_{m+k},\;\;{\rm if} \;\;m+k\ne
0;\;\;f(m,0)g(x_{-k})=c_kk(\beta-k)v_0.$$ Hence $g(x_m)=c_kv_{m+k}$,
when $m\neq-k$ and $g(x_{-k})=c_k(\beta-k)v_0$. In this case,
$\beta-k\ne 0$ since $g$ is an isomorphism. Since
$g(x_mx_n)=x_mg(x_n)$, we have
$$f(m,n)g(x_{m+n})=c_kx_mv_{n+k}\;\;{\rm if}\;\;n+k\ne 0,$$$$
{f(m,-k)}g(x_{m-k})=c_k(\beta-k)x_mv_0=0.$$

\noindent {\it Case (I-1).} $m+n+k \neq 0$. If $n\neq -k$, we have
$$f(m,n)c_k v_{m+n+k}=c_k x_m v_{n+k}=(n+k) c_k v_{m+n+k}.$$
Thus $f(m,n)=n+k$. If $n=-k$, we have $f(m,-k)=0$. Therefore we know
that
$$f(m,n)=n+k, \mbox{ for all } m,n \in \mathbb{Z} \mbox{ with } m+n+k \neq 0 .$$

\noindent {\it Case (I-2).} $m+n+k = 0$. If $n\neq -k$, we have
$$f(-k-n,n)c_k (\beta-k) v_0=c_k x_m v_{n+k}=(n+k)(\beta-n-k) c_k
{v_0.}$$ Then $f(-k-n,n)= \frac{(n+k)(\beta-n-k)}{ \beta-k}$. If
$n=-k$, we have $f(0,-k)=0$. Therefore we know that
$$f(m,n)=  \frac{(n+k)(\beta-n-k)}{\beta-k}, \mbox{ for
all } m,n \in \mathbb{Z} \mbox{ with }m+n+k = 0.$$

{\it Case (II).} $k=0$. Then $f(0,m)=m$. Since
$$x_0g(x_m)=g(x_0x_m)=f(0,m)g(x_m)=mg(x_m),$$
we know that $g(x_m)=d_mv_m$, where $d_m\neq0$. In particular,
$d_0=c_0$. Since
$$(n-m)d_{m+n}v_{m+n}=g([x_m,x_n])=x_mg(x_n)-x_ng(x_m)=d_nx_mv_n-d_mx_nv_m,$$
we have
$$(n-m)d_{m+n}=nd_n-md_m,\;\;{\rm if}\;\;   m+n \neq 0;\;\;
(\beta-n)d_n+(\beta+n)d_{-n}=2d_0.$$ Let $m=2$ and $n=-1$, we have
$d_{-1}=3d_1-2d_2$. Let $m=-2$ and $n=1$, we have
$d_{-2}=4d_1-3d_2$. On the other hand, let $n=1$, we have
$(\beta-1)d_1+(\beta+1)(3d_1-2d_2)=2d_0$. Let $n=2$, we have
$(\beta-2)d_2+(\beta+2)(4d_1-3d_2)=2d_0$. Therefore, we know that
$$d_0=\beta d_1,\;\;d_1=d_2=d_{-1}=d_{-2}.$$
By induction, we have $d_n=d_1$ for all $n \neq 0$ and $d_0=\beta
d_1$, where $\beta\neq0$ since $d_0\ne 0$. Because
$$f(m,n)g(x_{m+n})=g(x_mx_n)=x_mg(x_n),$$
with a similar discussion as in Case (I), we know that
$$
f(m,n)=\left\{
\begin{array}
{c@{ \quad }l}
n, & \mbox{if } m+n \neq 0; \\
 \frac{n(\beta -n)}{\beta}, & \mbox{if } m+n= 0, \mbox{
where } \beta\neq0.
\end{array} \right.$$
Obviously, it is just the case $k=0$ in equations (3.6) and (3.7).

Conversely, it is easy to know that we define a compatible
left-symmetric algebra structure on the Witt algebra $W$ by
equations (3.6) and (3.7).\hfill$\Box$

{\bf Lemma 3.9.} Let $T$ be an automorphism of the Witt algebra $W$.
Then $T(x_0)=\pm x_0$.

{\bf Proof.} Suppose $T(x_0)=\sum_ic_ix_i$ and $T(x_m)=\sum_jd_jx_j$
for a finite number of nonzero $c_i,d_j\in\mathbb{C}$. Therefore we
have
$$m\sum_jd_jx_j=T([x_0,x_m])=[T(x_0),T(x_m)]={\sum_{i+j=l}(j-i)c_id_jx_l}.$$
If there exists $i>0$ such that $c_i\neq0$, then let $n$ and $k$ be
the maximal numbers in the sets $\{i\mid c_i\neq0\}$ and $\{j\mid
d_j\neq0\}$ respectively. Therefore, by comparing the coefficient of
$x_{n+k}$ in both hand sides of the above equation, we have
$0=(k-n)c_nd_k$. So $k=n$. However, since $m$ is arbitrary, we know
that Im$T\subset \bigoplus_{i\leq n}\mathbb{C}x_i$, which is
contradictory  to the assumption that $T$ is an automorphism of $W$.
By a similar discussion, we can prove that there does not exist
$i<0$ such that $c_i\neq0$.

Therefore $T(x_0)=c_0x_0$. Since $m\sum_jd_jx_j=\sum_jjc_0d_jx_j$,
we know that there exists $k$ such that $d_k\neq0$ and  $m=kc_0$ for
any $m\neq0$. Hence $T(x_m)=d_kx_k$.  In particular, there exists
$m_1\in\mathbb{Z}$, such that $T(x_{m_1})=d_1x_1$. So $c_0=m_1$.
Moreover, since $c_0|m$ for any $m\in\mathbb{Z}$, we have $c_0=\pm
1$. \hfill $\Box$

{\bf Proposition 3.10.} Let $(V_1,\cdot)$ and $(V_2,*)$ be two
compatible left-symmetric algebras structures on the Witt algebra
$W$ defined by
$$x_m\cdot x_n=f_1(m,n)x_{m+n}, \quad  \mbox{ and }
\quad    x_m*x_n=f_2(m,n)x_{m+n}$$  respectively. If $V_1$ is
isomorphic to $V_2$ as left-symmetric algebras, then
$$f_1(m,n)=f_2(m,n),\;\;{\rm or}\;\;f_1(m,n)=-f_2(-m,-n).$$

{\bf Proof.} Let $T$ be an isomorphism from $V_1$ into $V_2$, that
is, $T(x\cdot y)=T(x)*T(y)$ for all $x,y\in V_1$. Obviously, $T$ is
an automorphism of $W$. By Lemma 3.9, we have $T(x_0)=\pm x_0$.

If  $T(x_0)=x_0$, then there exists $a_m\neq0$ such that $T(x_m)=a_m
x_m$ from the proof of Lemma 3.9. Hence it is obvious that $a_m
a_n=a_{m+n}$ and $a_0=1$. Since
$$f_1(m,n) a_{m+n}x_{m+n}=T(x_m\cdot x_n)=
T(x_m)*T(x_n)=f_2(m,n)a_m a_nx_{m+n},$$ we have $f_1(m,n)=f_2(m,n)$.

If  $T(x_0)=-x_0$, then there exists $b_m\neq0$ such that
$T(x_m)=b_m x_{-m}$ from the proof of Lemma 3.9. Hence it is easy to
know that $b_m b_n=-b_{m+n}$ and $b_0=-1$. Since
$$f_1(m,n)b_{m+n}x_{-(m+n)}=T( x_m\cdot x_n)=
T(x_m)*T(x_n)=f_2(-m,-n) b_m b_n x_{-(m+n)},$$ we have
$f_1(m,n)=-f_2(-m,-n)$. \hfill $\Box$

{\bf Theorem 3.11.} Any graded compatible left-symmetric algebra
structure on the Witt algebra $W$ satisfying equation (3.1) is
isomorphic to one of the following algebras:
$$V_{\alpha,\epsilon},\;\;\alpha,\epsilon \in
\mathbb{C}\;\; {\rm satisfying}\;\; \epsilon=0\;\; {\rm or}\;\;
\epsilon^{-1}\notin\mathbb{Z};$$$$V^{\beta,k},\;\;\beta\in
\mathbb{C}\;\;{\rm and}\;\; k \in \mathbb{Z}\;\;{\rm satisfying}\;\;
\beta \ne k.$$ Moreover, the isomorphisms between them are exactly
given as follows,
$$V_{\alpha,\epsilon}\cong V_{-\alpha,-\epsilon},\ \
V_{\alpha,0}\cong V_{\alpha,1/\alpha} \mbox{ with
}\alpha\notin\mathbb{Z}\ \mbox{  and  }\ V^{\beta,0}\cong
V^{-\beta,0}.$$

{\bf Proof.} The first half part immediately follows from Corollary
3.4 and Theorems 3.5-3.8. For the second half part, we assume that
$V_1\cong V_2$, where $V_1, V_2$ are the compatible left-symmetric
algebra structures on $W$ defined by equation (3.1) with $f_1(m,n)$
and $f_2(m,n)$ respectively. Then by Proposition 3.10, we have
$f_1(m,n)=f_2(m,n)$ or $f_1(m,n)=-f_2(-m,-n)$. Obviously, in the
former case, $V_1=V_2$ and in the latter case, $V_1\cong V_2$ by the
linear isomorphism $x_n\rightarrow -x_{-n}$.

 {\it Case (I).} $V_1=V_{\alpha_1,\epsilon_1}$ and
$V_2=V_{\alpha_2,\epsilon_2}$.

If $f_1(m,n)=f_2(m,n)$, then $f_1(0,0)=\alpha_1=f_2(0,0)=\alpha_2$.
Set $\alpha_1=\alpha_2=\alpha$. Since
$$f_1(m,n)=\frac{(\alpha+n+\alpha\epsilon_1 m)(1+\epsilon_1
n)}{1+\epsilon_1(m+n)}=\frac{(\alpha+n+\alpha\epsilon_2
m)(1+\epsilon_2 n)}{1+\epsilon_2(m+n)}=f_2(m,n),$$ we have
$$(\epsilon_1-\epsilon_2)[\alpha(\epsilon_1+\epsilon_2)-1+\alpha\epsilon_1\epsilon_2(m+n)]=0.$$
If $\alpha=0$, then $\epsilon_1=\epsilon_2$. In this case,
$V_1=V_2={V_{0,\epsilon}}$. If $\alpha\ne 0$, then
$\epsilon_1=\epsilon_2$ or $\epsilon_1=0$,
$\epsilon_2=\frac{1}{\alpha}$ or $\epsilon_2=0$,
$\epsilon_1=\frac{1}{\alpha}$. Therefore, we have
$V_1=V_2=V_{\alpha,\epsilon_1},$ or $$ V_1=V_{\alpha,
0},\;\;V_2=V_{\alpha,1/\alpha},\alpha\notin\mathbb{Z};\;\;{\rm
or}\;\;V_1=V_{\alpha,1/\alpha}, V_2=V_{\alpha,
0},\;\;\alpha\notin\mathbb{Z}.$$

If $f_1(m,n)=-f_2(-m,-n)$, then
$f_1(0,0)=\alpha_1=-f_2(0,0)=-\alpha_2$. Set
$\alpha_1=-\alpha_2=\alpha$. Since
$$f_1(m,n)=\frac{(\alpha+n+\alpha\epsilon_1 m)(1+\epsilon_1
n)}{1+\epsilon_1(m+n)}=\frac{(\alpha+n-\alpha\epsilon_2
m)(1-\epsilon_2 n)}{1-\epsilon_2(m+n)}=-f_2(-m,-n),$$ we have
$$(\epsilon_1+\epsilon_2)[\alpha(\epsilon_1-\epsilon_2)-1-\alpha\epsilon_1\epsilon_2(m+n)]=0.$$
If $\alpha=0$, then $\epsilon_1=-\epsilon_2$. In this case,
$V_1=V_{0,\epsilon}\cong V_2=V_{0,-\epsilon}$. If $\alpha\ne 0$,
then $\epsilon_1=-\epsilon_2$ or $\epsilon_1=0$,
$\epsilon_2=-\frac{1}{\alpha}$ or $\epsilon_2=0$,
$\epsilon_1=\frac{1}{\alpha}$. Therefore, we have
$$V_1=V_{\alpha, \epsilon},\;{V_2=V_{-\alpha, -\epsilon}};\;{\rm or}\;\; V_1=V_{\alpha,
0},\;V_2=V_{-\alpha,-1/\alpha},\alpha\notin\mathbb{Z};\;$$ $${\rm
or}\;V_1=V_{\alpha,1/\alpha}, V_2=V_{-\alpha,
0},\;\alpha\notin\mathbb{Z}.$$

{\it Case (II).} $V_1=V^{\beta_1,k_1}$ and $V_2=V^{\beta_2,k_2}$.

If $f_1(m,n)=f_2(m,n)$, then $f_1(0,0)=k_1=f_2(0,0)=k_2$. Set
$k_1=k_2=k$. Since
$$f_1(-n-k,n)=\frac{(n+k)(\beta_1-n-k)}{\beta_1-k}=\frac{(n+k)(\beta_2-n-k)}{\beta_2-k}=f_2(-n-k,n),$$
we know that $\beta_1=\beta_2$.

If $f_1(m,n)=-f_2(-m,-n)$, then $f_1(0,0)=k_1=-f_2(0,0)=-k_2$. Set
$k_1=-k_2=k$. If $k\ne 0$, then let $n\neq 0,-k$. Since
$f_1(-n-k,n)=-f_2(n+k,-n)$, we have
$$ \frac{(n+k)(\beta_1-n-k)}{\beta_1-k}=-(-n-k).$$
So $\beta_1-n-k=\beta_1-k$, which is a contradiction. Therefore
$k=0$. Since
$$f_1(-n,n)= \frac{n(\beta_1-n)}{\beta_1}=-\frac{-n(\beta_2+n)}{\beta_2}=-f_2(n,-n),$$
we have $\beta_1=-\beta_2$. On the other hand, when
$\beta_1=-\beta_2=\beta$ and $k_1=k_2=0$, it is easy to know that
$f_1(m,n)=-f_2(-m,-n)$ and then $V_1=V^{\beta,0}\cong
V_2=V^{-\beta,0}$.

{\it Case (III).} $V_1=V_{\alpha,\epsilon}$ and $V_2=V^{\beta,k}$.

If $f_1(m,n)=f_2(m,n)$, then $f_1(0,0)=\alpha=f_2(0,0)=k\in
\mathbb{Z}$. Obviously there exist $m,n$ such that $m\ne 0$, $n\ne
0$ and $m+n+k\neq0$. Since
$$f_1(m,n)=\frac{(k+n+k\epsilon m)(1+\epsilon
n)}{1+\epsilon(m+n)}=n+k=f_2(m,n),$$ we have
$$\epsilon mn(k\epsilon-1)=0.$$
Since $\epsilon^{-1}\notin\mathbb{Z}$, we know $k\epsilon-1\neq0$.
Then $\epsilon=0$. Let $n\neq0,-k$. Since
$$f_1(-n-k,n)=k+n=\frac{(k+n)(\beta-k-n)}{\beta-k}=f_2(-n-k,n),$$
we have $\beta-k-n=\beta-k$, which is a contradiction.

If $f_1(m,n)=-f_2(-m,-n)$, then we still have $f_1(m,n)=f_2(m,n)$ by taking
$V_1=V_{-\alpha,-\epsilon}$ and $V_2=V^{\beta,k}$. Since
$V_{\alpha,\epsilon}\cong V_{-\alpha,-\epsilon}$, we know that
$V_{\alpha,\epsilon}$ is not isomorphic to  $V^{\beta,k}$ in this subcase.

So there does not exist an isomorphism between $V_{\alpha,\epsilon}$
and $V^{\beta,k}$. \hfill$\Box$

{\bf Example 3.12.} In \cite{Cha}, the notion of pre-Lie algebra was
used. Chapoton classified the simple graded left-symmetric algebras
of growth one satisfying the following two conditions:

(1) The underlying graded vector space is $E=\oplus_{i\in \mathbb
Z}{\mathbb C} e_i$;

(2) The product is given by {$e_i\circ
e_j=f(i)g(j)e_{i+j}$}.

\noindent Such a left-symmetric algebra is isomorphic
either to the algebra $A_a$ defined by
$$e_i\circ e_j=(1+aj)e_{i+j},\;\;a\in {\mathbb C},$$
or to the algebra $B_b$ defined by
$$e_i\circ e_j=\frac{j}{1+bi}e_{i+j},\;\; b=0\;{\rm or}\;b^{-1}\notin \mathbb Z.$$
Obviously, $A_0$ is a commutative associative algebra and $B_0\cong
V_{0,0}$. Moreover, when $a\ne 0$, $A_a\cong V_{\frac{1}{a},0}$ by a
linear transformation $e_i\rightarrow \frac{1}{a} e_i$. When $b\ne
0, b^{-1}\notin \mathbb Z$, $B_b\cong V_{0, b}$ by a linear
transformation $e_i\rightarrow (1+bi)e_i$. Therefore, the only
isomorphisms among $A_a$ and $B_b$ are $A_a\cong A_{-a}$ and
$B_b\cong B_{-b}$ due to Theorem 3.11, which were given in [Cha],
too. Hence, except $A_0$, we have obtained all left-symmetric
algebras in [Cha].

In particular, $V_{0,b}$ for $b=0$ or $b^{-1}\notin \mathbb Z$ was also given in [Ku2] as a special
case of the left-symmetric algebras satisfying equation (3.1).  \hfill $\Box$

{\bf Example 3.13.} Recall that a Novikov algebra $A$ is a
left-symmetric algebra satisfying $(xy)z=(xz)y$ for any $x,y,z\in
A$. Novikov algebras have been introduced in connection with
Hamiltonian operators in the formal variational calculus (\cite{GD})
and Poisson brackets of hydrodynamic type (\cite{BN}). It is easy to
see that a compatible Novikov algebra structure on the Witt algebra
satisfying equation (3.1) must be isomorphic to one of
$N_\alpha=V_{\alpha,0}$, that is,
$$x_mx_n=(\alpha+n)x_{m+n}, \quad \alpha\in\mathbb{C}.$$
Moreover, the isomorphisms between them are exactly given by
$N_\alpha\cong N_{-\alpha}$. It is just the case (ii) appearing in
\cite{O}.    \hfill $\Box$

\section{Compatible left-symmetric algebra structures on the Virasoro algebra}

In this section, we consider the central extensions of the
left-symmetric algebras obtained in Section 3 whose commutator is
the Virasoro algebra $\mathscr{V}$.

Let $(A,\cdot)$ be a left-symmetric algebra and $\omega:A\times
A\rightarrow \mathbb{C}$ be a bilinear form. It defines a
multiplication on the space $\widehat{A}=A\oplus \mathbb{C}\theta$,
by the rule
$$(x+a\theta)\ast(y+b\theta)=x\cdot y +\omega(x,y)\theta,\quad\quad x,y\in A,\
a,b\in \mathbb{C}.\eqno{(4.1)}$$ Obviously, $\widehat{A}$ is a
left-symmetric algebra if and only if
$$\omega(x\cdot y,z)-\omega(x,y\cdot z)=\omega(y\cdot
x,z)-\omega(y,x\cdot z).\eqno(4.2)$$ $\widehat A$ is called a
central extension of $A$. Moreover, by construction, the bilinear
form
$$\Omega(x,y)=\omega(x,y)-\omega(y,x),\quad\mbox{ where }
 x,y\in {\mathcal G}(A)\,,\eqno{(4.3)}$$
defines a central extension of its sub-adjacent Lie algebra ${\mathcal G}(A)$.

Let $V$ be a compatible left-symmetric algebra structure on the Witt
algebra $W$ satisfying equation (3.1). Since the Virasoro algebra
$\mathscr{V}$ is a central extension of the Witt algebra $W$, it is
natural to consider the central extension $\widehat V=V\oplus
\mathbb C\theta$ of $V$ such that $\widehat V$ is a compatible
left-symmetric algebra structure on the Virasoro algebra
$\mathscr{V}$ while $\theta$ being the annihilator of $\mathscr{V}$,
that is, the products are given by
$${\theta\theta=}x_m\theta=\theta
x_m=0 \mbox{ and } x_mx_n=f(m,n)x_{m+n}+\omega(x_m,x_n)\theta \eqno
(4.4)$$ where $f(m,n)$ satisfies equations (3.3)-(3.4).

We denote $\omega(x_m,x_n)$ by $\omega (m,n)$ for convenience. Then
by equations (2.8), (4.2) and (4.3), we have the following equations
$$\omega(m,n)-\omega(n,m)=\frac{1}{12}(n^3-n)\delta_{m+n,0},\eqno(4.5)$$
$$(n-m)\omega(m+n,l)=\omega(m,n+l)f(n,l)-\omega(n,m+l)f(m,l).\eqno(4.6)$$

{\bf Theorem 4.1.} When $\alpha\neq0$ or $\epsilon=0$, there does
not exist a central extension of $V_{\alpha,\epsilon}$ satisfying
equations (4.5)-(4.6). There is exactly one central extension of
$V_{0,\epsilon}$ with $\epsilon\neq0$ and $\epsilon^{-1} \notin
\mathbb{Z}$ satisfying equations (4.5)-(4.6), which is given by:
$$\omega(x_m,x_n)=\omega(m,n)=\frac{1}{24}(n^3-n-(\epsilon-\epsilon^{-1})n^2)\delta_{m+n,0}.\eqno (4.7)$$

{\bf Proof.} Let $m=-n\neq0,\ l=0$ in equation (4.6), we have
$$\omega(0,0)=\frac{1}{24}\alpha(n^2-1).$$
Since $\omega(0,0)$ does not depend on $n$, we know that $\alpha=0$ and $\omega(0,0)=0$. Notice that
$V_{0,\epsilon}$ is given by (see [Ku2], too)
$$f(m,n)= \frac{n(1+\epsilon n)}{1+\epsilon(m+n)}.$$
Let $n=l=0$ in equations (4.5) and (4.6), then we get
$$\omega(m,0)=\omega(0,m)=0.$$
Let $m=0$ in equation (4.6), thus we know
$$n\omega(n,l)=\omega(0,n+l)f(n,l)-\omega(n,l)f(0,l)=-l\omega(n,l),$$
that is, $(n+l)\omega(n,l)=0$. So we can assume that
$$\omega(n,l)=\varphi(n)\delta_{n+l,0} \mbox{ for some map }\
\varphi:\mathbb{Z}\rightarrow\mathbb{C}.$$ Let $m+n=0$ in equation
(4.5), thereby we have
$$\varphi(n)-\varphi(-n)=\frac{1}{12}(-n^3+n).$$
Let $m+n+l=0$ in equation (4.6), then
$$(n-m)\varphi(m+n)=\varphi(m)f(n,-m-n)-\varphi(n)f(m,-m-n),$$
which gives
$$(n-m)\varphi(m+n)=\frac{(-m-n)(1-\epsilon(m+n))}{1-\epsilon
m}\varphi(m)-\frac{(-m-n)(1-\epsilon(m+n))}{1-\epsilon
n}\varphi(n).$$ Set $\psi(m)= \frac{\varphi(m)}{m(1-\epsilon m)}$.
Then we know that
$$\left\{\begin{split}&(1-\epsilon n)\psi(n)+(1+\epsilon n)\psi(-n)=\frac{1}{12}(-n^2+1),\\
&(n-m)\psi(m+n)=-m\psi(m)+n\psi(n). \end{split}\right.\eqno (***)$$
Let $m=2, n=-1$ in equation ($***$), so we have
$$(1+\epsilon)\psi(-1)+(1-\epsilon)\psi(1)=0,\;\;-3\psi(1)=-2\psi(2)-\psi(-1).$$
So $\psi(2)=\frac{\epsilon+2}{\epsilon+1}\psi(1)$.  Let $m=-2, n=1$
in the first part of equation ($***$), hence we get
$$3\psi(-1)=2\psi(-2)+\psi(1).$$
So $\psi(-2)=\frac{\epsilon-2}{\epsilon+1}\psi(1)$. Notice that (by the first part of equation ($***$))
$$(1-2\epsilon)\psi(2)+(1+2\epsilon)
\psi(-2)=- \frac{1}{4}.$$ If $\epsilon=0$, then
$$\psi(2)+\psi(-2)=2\psi(1)-2\psi(1)=0= - \frac{1}{4},$$
which is a contradiction. So we suppose that $\epsilon\neq0$. Then
$\psi(1)= \frac{ 1}{24}(1+\epsilon^{-1})$. Let $n=1$ in the second
part of equation ($***$), we have
$$(1-m)\psi(m+1)=-m\psi(m)+\psi(1).$$
Then
$$(m-1)(\psi(m+1)-\psi(1))=m(\psi(m)-\psi(1)).$$
Without losing generality, we can assume that $m\geq 2$ (a similar
discussion is for $m\leq 0$). Therefore, $$\psi(m+1)-\psi(1)=
\frac{m}{m-1}\cdot\frac{m-1}{m-2}\cdots\frac{2}{1}(\psi(2)-\psi(1))
=m(\psi(2)-\psi(1)){=\frac{m}{1+\epsilon}\psi(1)}.$$ Hence
$$\psi(m)= \frac{m-1}{1+\epsilon}\psi(1)+\psi(1)=
\frac{1}{24}(1+\epsilon^{-1}m).$$ So
$$\varphi(m)= \frac{1}{24}m(1-\epsilon
m)(1+\epsilon^{-1}m),$$ and then
$$\omega(m,n)=\varphi(-n)\delta_{m+n,0}=\frac{1}{24}(n^3-n-(\epsilon-\epsilon^{-1})n^2)\delta_{m+n,0}.$$
Moreover, it is easy to know that $\omega(m,n)$ satisfies equations
(4.5)-(4.6). \hfill $\Box$

{\bf Theorem 4.2.} There does not exist a  central extension of
$V^{\beta,\epsilon}$ satisfying (4.5)-(4.6).

{\bf Proof.}  Firstly, let $m=-n\neq 0$ and $l=0$ in equation (4.6),
we can get
$$\omega(0,0)=\frac{1}{24}k(n^2-1).$$
Since $\omega (0,0)$ does not depend on $n$, we know that $k=0$ and
$\omega (0,0)=0$. Let $n=l=0$ in equations (4.5) and (4.6), we have
$$\omega(m,0)=\omega(0,m)=0.$$
On the other hand, let $m=0$ in equation (4.6), then we have
$$n\omega(n,l)=\omega(0,n+l)f(n,l)-\omega(n,l)f(0,l)=-l\omega(n,l),$$
that is, $(n+l)\omega(n,l)=0$. So we can assume that
$$\omega(n,l)=\varphi(n)\delta_{n+l,0} \mbox{ for some map }\
\varphi:\mathbb{Z}\rightarrow\mathbb{C}.$$ Finally, let $m+n=0$ in
equation (4.5), thus we get
$$\varphi(n)-\varphi(-n)=\frac{1}{12}(-n^3+n).$$
So $\varphi(1)-\varphi(-1)=0$. Let $m,n\neq0,m+n+l=0$ in equation
(4.6), thereby we know
$$(n-m)\varphi(m+n)=(m+n)(\varphi(n)-\varphi(m)).$$
Let $m=2, n=-1$ and $m=-2,n=1$ in the above equation respectively, we know that
$$\varphi(2)= \varphi(-2)=4\varphi(1).$$
Hence $\varphi(2)-\varphi(-2)=0$. However, on the other hand,
$$\varphi(2)-\varphi(-2)=\frac{1}{12}(-2^3+2)=-\frac{1}{2},$$
which is a contradiction.\hfill$\Box$

{\bf Corollary 4.3.} Any compatible left-symmetric algebra structure
on the Virasoro algebra $\mathscr V$ satisfying equation (4.4) is
isomorphic to one of the algebras given by the multiplication
$$x_mx_n=\frac{n(1+\epsilon
n)}{1+\epsilon(m+n)}x_{m+n}+\frac{\theta}{24}(n^3-n-(\epsilon-\epsilon^{-1})n^2)\delta_{m+n,0},\eqno
(4.8)$$ where $\theta$ is an annihilator and $\rm{Re} \epsilon>0,
\epsilon^{-1}\notin\mathbb{Z}\ \mbox{ or }\ \rm{Re}
\epsilon=0,\rm{Im} \epsilon>0.$

{\bf Remark 4.4.} The equation (4.8) was also obtained in [Ku2] as a
central extension of the left-symmetric algebra $V_{0,\epsilon}$ for
$\epsilon\ne 0, \epsilon^{-1}\notin \mathbb Z$, which was given only
as a special case of the left-symmetric algebras satisfying equation
(3.1). It is interesting to see that in our discussion there does
not exist a central extension satisfying equations (4.5)-(4.6) in
other cases.

\section*{Acknowledgements} The authors thank Professor S.B. Tan for
valuable discussions. This work was supported by the National
Natural Science Foundation of China (10571091,10621101), NKBRPC
(2006CB805905), Program for New Century Excellent Talents in
University.

\end{document}